\documentclass[fleqn]{mat01}
\usepackage{times,mathtimy,amssymb,latexsym}
\begin{document}



\setcounter{page}{233}

\firstpage{233}

\markboth{Leszek Gasi\'{n}ski, Dumitru Motreanu and Nikolaos S
Papageorgiou}{Multiplicity of nontrivial solutions}

\title{Multiplicity of nontrivial solutions for elliptic
equations with nonsmooth potential and resonance at higher
eigenvalues}

\author{LESZEK GASI\'NSKI$^{1}$, DUMITRU MOTREANU$^{2}$ and~NIKOLAOS~S~PAPAGEORGIOU$^{3}$}

\address{$^{1}$Institute of Computer Science, Jagiellonian University,
ul.~Nawojki 11, 30072~Cracow, Poland\\
\noindent $^{2}$D\'epartement de Mathematiques, Universit\'e de
Perpignan, 66860~Perpignan, France\\
\noindent $^{3}$Department of Mathematics, National Technical
University, Zografou Campus, Athens~15780, Greece\\
\noindent E-mail: npapg@math.ntua.gr}

\volume{116}

\mon{May}

\parts{2}

\pubyear{2006}

\Date{MS received 15 March 2004; revised 22 February 2005}

\begin{abstract}
We consider a semilinear elliptic equation with a nonsmooth,
locally \hbox{Lipschitz} potential function (hemivariational
inequality). Our hypotheses permit double resonance at infinity
and at zero (double-double resonance situation). Our approach is
based on the nonsmooth critical point theory for locally Lipschitz
functionals and uses an abstract multiplicity result under local
linking and an extension of the Castro--Lazer--Thews reduction
method to a nonsmooth setting, which we develop here using tools
from nonsmooth analysis.
\end{abstract}


\maketitle

\section{Introduction} \label{Intr}

Let $Z\subseteq\mathbb{R}^N$ be a bounded domain with a $C^2$
boundary $\Gamma$. We study the following resonant semilinear
elliptic differential equation with a nonsmooth potential
(hemivariational inequality):
\begin{equation*}
\begin{cases}
-\Delta x(z)-\lambda_k x(z) \in \partial j(z,x(z)),
\ \ \ \hbox{for a.a.}\ z\in Z   \\[.4pc]
x|_{\Gamma}=0.
\end{cases}
\tag{\rm HVI}
\end{equation*}

Here $k\geq 1$ is a fixed integer, $\{\lambda_n\}_{n\geq 1}$ is
the increasing sequence of distinct eigenvalues of the negative
Laplacian with Dirichlet boundary condition (i.e. of
$(-\Delta,H^1_0(Z))$), $j(z,\zeta)$ is a locally Lipschitz in the
$\zeta$-variable integrand (in general it can be nonsmooth) and
$\partial j(z,\zeta)$ is the Clarke subdifferential with respect
to the $\zeta$-variable.

For problem (HVI), we prove a multiplicity result, using a recent
abstract theorem on the existence of multiple nontrivial critical
points for a nonsmooth locally Lipschitz functional, proved by
Kandilakis, Kourogenis and Papageorgiou~\cite{KKP}. Our approach
is variational and is based on the nonsmooth critical point theory
(see \cite{Ch} and \cite{KP}). In particular, we develop and use a
nonsmooth variant of the so-called `reduction method'. This method
was first introduced for smooth problems by Castro and
Lazer~\cite{CL} and Thews~\cite{T} (see also \cite{CC,CCN}). Our
hypotheses allow the nonsmooth potential to interact
asymptotically at $\pm\infty$ with two consecutive eigenvalues of
higher order. Berestycki and de~Figueiredo~\cite{BF} were the
first to consider such problems with resonance between $\lambda_1$
and $\lambda_2$ and they coined the term `double resonance
problems'. In their analysis, the use of the interval
$[\lambda_1,\lambda_2]$ is crucial, since they exploit heavily the
fact that the principal eigenfunction $u_1$ is strictly positive
and ${\partial u_1}/{\partial n}<0$ with $n$ being the outward
unit normal on the boundary (this is a consequence of the strong
maximum principle). It is well-known that in higher parts of the
spectrum this is no longer true. Recall that the principal
eigenfunction $u_1$ is the only one with constant sign. So in
higher parts of the spectrum the analysis is more delicate.
Recently for smooth problems, the issue was investigated by
Cac~\cite{Ca}, Hirano and Nishimura~\cite{HN}, Robinson~\cite{Ro},
Costa-Silva~\cite{CS}, Landesman, Robinson and Rumbos~\cite{LRR},
Iannacci and Nkashama~\cite{IN}, Tang and Wu~\cite{TWa}, Su and
Tang~\cite{ST} and Su~\cite{S}.

For problems with nonsmooth potential (known in the literature as
hemivariational inequalities), equations resonant at higher
eigenvalues were investigated by Goeleven, Motreanu and
Panagiotopoulos~\cite{GMP} and Gasi\'nski and
Papageorgiou~\cite{GPc}. However, they did not allow for the
situation of double resonance.

We should mention that hemivariational inequalities arise in
physical problems, when one wants to consider more realistic
models with a nonsmooth and nonconvex energy functionals. For
concrete applications we refer to the book of Naniewicz and
\hbox{Panagiotopoulos}~\cite{NP}. For the mathematical theory of
hemivariational inequalities we refer to the work of Gasi\'nski
and Papageorgiou~\cite{GPa,GPb}, Motreanu and
Panagiotopoulos~\cite{MPa,MPb}, Niculescu and Radulescu~\cite{NR},
Radulescu~\cite{R}, Radulescu and Panagiotopoulos~\cite{RP} and
the references therein.

\section{Mathematical background}\label{prelim}

As we have already mentioned, our approach is based on the theory
of the nonsmooth critical point theory for locally Lipschitz
functionals.

Let $X$ be a Banach space and $X^*$ its topological dual. By
$\|\cdot\|_X$ we denote the norm of $X$ and by
$\langle\cdot,\cdot\rangle_X$ the duality pairing for the pair
$(X,X^*)$.

We will be dealing with locally Lipschitz functions
$\varphi\hbox{\rm :}\ X\rightarrow \mathbb{R}.$

Recall that a continuous convex function is locally Lipschitz. For
a locally Lipschitz function $\varphi\hbox{\rm :}\ X\rightarrow
\mathbb{R}$, we introduce the {\it generalized directional
derivative} of $\varphi$ at $ x\in X$ in the direction  $h\in X$,
defined by
\begin{equation*}
\varphi^0(x;h)\stackrel{d f}{=}\limsup\limits_{\substack
{x' \,\rightarrow \,x\\[.4pc]
t \,\searrow \,0}} \frac{\varphi(x'+th)-\varphi(x')}{t}
\end{equation*}
(see \cite{Cl}). It is easy to check that the function $X\ni
h\rightarrow\varphi^0(x;h)\in \mathbb{R}$ is sublinear,
continuous, and so by the Hahn--Banach theorem,
$\varphi^0(x;\cdot)$ is the support function of a nonempty, convex
and $w^*$-compact set $\partial\varphi(x)$, defined by
\begin{equation*}
\partial\varphi(x)\stackrel{d f}{=}\{ x^*\in X^*\hbox{\rm :}\
\langle x^*,h\rangle_X\leq \varphi^0(x;h) \quad \hbox{for all}\
h\in X\}.
\end{equation*}
The multifunction $\partial\varphi\hbox{\rm :}\ X\rightarrow
2^{X^*}\!\setminus\!\{\emptyset\}$ is known as the {\it
generalized} (or {\it Clarke}) {\it subdifferential} of $\varphi$.
If $\varphi\in C^1(X),$ then $\partial \varphi
(x)=\{\varphi'(x)\}$ and if $\varphi$ is convex, then $\partial
\varphi(x)$ coincides with the convex subdifferential. The
generalized subdifferential has a rich calculus which can be found
in \cite{Cl}.

A point $x\in X$ is a critical point of the locally Lipschitz
function $\varphi$, if $0\in\partial\varphi(x)$. If $x\in X$ is a
critical point, the value $c=\varphi(x)$ is a critical value of
$\varphi$. It is easy to check that if $x\in X$ is a local
extremum of $\varphi$ (i.e. a local minimum or a local maximum),
then $0\in\partial\varphi(x)$ (i.e. $x\in X$ is a critical point).

We will use the following compactness-type condition:\vspace{.8pc}

\noindent A locally Lipschitz function $\varphi\hbox{\rm :}\
X\rightarrow \mathbb{R}$ satisfies the {\it nonsmooth
Palais--Smale condition} ({\it nonsmooth {\rm PS}-condition} for
short) if any sequence $\{ x_n\}_{n\geq 1}\subseteq X$ such that
$\{\varphi(x_n)\}_{n\geq 1}$ is bounded and
\begin{equation*}
m_{\varphi}(x_n)\stackrel{d f}{=} \min \{ \|x^*\|_{X^*}\hbox{\rm
:}\ x^*\in \partial\varphi(x_n)\} \longrightarrow 0\quad
\hbox{as}\ n\rightarrow +\infty,
\end{equation*}
has a strongly convergent subsequence.\vspace{.8pc}

Since for $\varphi\in C^1(X)$ we have
$\partial\varphi(x)=\{\varphi'(x)\}$, we see that the above
definition is an extension of the smooth PS-condition. It was
shown by Cerami~\cite{Ce} and Bartolo, Benci and
Fortunato~\cite{BBF} that a slightly more general condition in the
smooth setting, suffices to prove the main minimax principles. In
the present nonsmooth setting this condition has the following
form:\vspace{.8pc}

\noindent A locally Lipschitz function $\varphi\hbox{\rm :}\
X\rightarrow \mathbb{R}$ satisfies the {\it nonsmooth Cerami
condition} ({\it non\-smooth {\rm C}-condition} for short) if any
sequence $\{ x_n\}_{n\geq 1}\subseteq X$ such that
$\{\varphi(x_n)\}_{n\geq 1}$ is bounded and
\begin{equation*}
(1+\|x_n\|_X) m_{\varphi}(x_n) \longrightarrow 0 \quad \hbox{as}\
n\rightarrow +\infty,
\end{equation*}
has a strongly convergent subsequence.\vspace{.8pc}

Recently Kandilakis, Kourogenis and Papageorgiou~\cite{KKP},
proved the following multiplicity result extending a corresponding
theorem of Brezis and Nirenberg~\cite{BN}.

\renewcommand\theequation{\thesection\arabic{equation}}

\begin{theorem}[\!] \label{thm1}
If $X$ is a reflexive Banach space{\rm ,} $X=Y\oplus V$ with $\dim
V<+\infty${\rm ,} $\varphi\hbox{\rm :}\ X\rightarrow\mathbb{R}$ is
locally Lipschitz{\rm ,} bounded below{\rm ,} satisfies the
nonsmooth {\rm C}-condition{\rm ,} $\inf_X \varphi<0$ and there
exists $\varrho>0$ such that
\begin{equation} \label{eql_543}
\left\{\begin{array}{ll}
\varphi(x)\leq 0 &\hbox{if}\ \ x\in V,\ \|x\|_X\leq\varrho\\[.4pc]
\varphi(x)\geq 0 &\hbox{if}\ \ x\in Y,\ \|x\|_X\leq\varrho
\end{array},\right.
\end{equation}
then $\varphi$ has at least two nontrivial critical points.
\end{theorem}

Condition~(\ref{eql_543}) implies $\varphi(0)=0$. We call
condition~(\ref{eql_543}) the local linking condition.

Recall that, if $\{\lambda_n\}_{n\geq 1}$ are the distinct
eigenvalues of $(-\Delta,H^1_0(Z))$, then
$\lambda_n\longrightarrow+\infty$ and $\lambda_1$ is positive,
simple and isolated. Also there is an orthonormal basis
$\{u_n\}_{n\geq 1}\subseteq H^1_0(Z)\cap C^{\infty}(Z)$ of
$L^2(Z)$, which are eigenfunctions corresponding to the
eigenvalues $\{\lambda_n\}_{n\geq 1}$, i.e.
\begin{equation*}
\left\{\begin{array}{l}
-\Delta u_n(z)=\lambda_nu_n(z) \quad \forall z\in Z,\\[.4pc]
u_n|_{\Gamma}=0,
\end{array}\right.
\end{equation*}
for $n\geq 1$. If the boundary $\Gamma$ of $Z$ is a $C^k$-manifold
(respectively a $C^{\infty}$-manifold) then $u_n\in C^k(\bar{Z})$
(respectively $u_n\in C^{\infty}(\bar{Z})$). The sequence
$\big\{\frac{1}{\sqrt{\lambda_n}}u_n\big\}_{n\geq 1}$ is an
orthonormal basis of $H^1_0(Z)$. For every integer $m\geq 1$, let
$E(\lambda_m)$ be the eigenspace corresponding to the eigenvalue
$\lambda_m$. We define
\begin{equation*}
\bar{H}_m \stackrel{d f}{=} \bigoplus_{i=1}^{m-1} E(\lambda_i)
\quad\hbox{and}\quad \hat{H}_m \stackrel{d f}{=}
\overline{\bigoplus_{i=m+1}^{\infty} E(\lambda_i)}.
\end{equation*}
We have the following orthogonal direct sum decomposition:
\begin{equation*}
H^1_0(Z) =  \bar{H}_m\oplus E(\lambda_m)\oplus\hat{H}_m.
\end{equation*}
The eigenspace $E(\lambda_m)\subseteq H^1_0(Z)\cap C^{\infty}(Z)$
has the unique continuation property, namely if $u\in
E(\lambda_m)$ is such that $u$ vanishes on a set of positive
measure, then $u(z)=0$ for all $z\in Z$.

If we set
\begin{equation*}
V_m \stackrel{d f}{=} \bar{H}_m\oplus E(\lambda_m)
\quad\hbox{and}\quad W_m \stackrel{d f}{=}
E(\lambda_m)\oplus\hat{H}_m,
\end{equation*}
then on these spaces we have variational characterizations of the
eigenvalues (Rayleigh quotients), which can be found in
\cite{jost}.

Let us recall two modes of convergence of sets and functions,
which will be used in the proof of our nonsmooth extension of the
Castro--Lazer--Thews reduction method. So let $(Y_1,\tau_1)$ and
$(Y_2,\tau_2)$ be two Hausdorff topological spaces ($\tau_1$ and
$\tau_2$ being the respective topologies). Also let
$\{G_n\}_{n\geq 1}$ be a sequence of nonempty subsets of
$Y_1\times Y_2$. We define:
\begin{align*}
(\tau_1\times\tau_2) - \liminf_{n\rightarrow+\infty}G_n
&\stackrel{df}{=} \Big\{ (u,v)\in Y_1\times Y_2\hbox{\rm :}\
u=\tau_1 - \lim_{n\rightarrow+\infty} u_n,\\[.4pc]
&\qquad\, v=\tau_2 - \lim_{n\rightarrow+\infty} v_n,
(u_n,v_n)\in G_n, n\geq 1 \Big\},\\[.4pc]
(\tau_1\times\tau_2) - \limsup_{n\rightarrow+\infty}G_n
&\stackrel{df}{=} \Big\{ (u,v)\in Y_1\times Y_2\hbox{\rm :}\
u=\tau_1 - \lim_{n\rightarrow+\infty} u_{n_k},\\[.4pc]
&\qquad\, v=\tau_2 - \lim_{n\rightarrow+\infty} v_{n_k},
(u_{n_k},v_{n_k})\in G_{n_k},\\[.4pc]
&\qquad\, n_i<n_{i+1}\ \hbox{for}\ i\geq 1 \Big\}.
\end{align*}
If
\begin{equation*}
G=(\tau_1\times\tau_2) - \liminf_{n\rightarrow+\infty}G_n=
(\tau_1\times\tau_2) - \limsup_{n\rightarrow+\infty}G_n,
\end{equation*}
then we say that the sequence $\{G_n\}_{n\geq 1}$
converges in the $(\tau_1\times\tau_2)$-sequential Kuratowski
sense to $G$ and denote it by
\begin{equation*}
G_n\stackrel{K_{\tau_1,\tau_2}}{\longrightarrow} G \quad\hbox{as}\
n\rightarrow+\infty.
\end{equation*}

Now, let $Y$ be a Banach space and $\{\varphi_n\}_{n\geq
1}\subseteq\bar{\mathbb{R}}^Y$ and $\varphi\in\bar{\mathbb{R}}^Y$
(with $\bar{\mathbb{R}}\stackrel{d
f}{=}\mathbb{R}\cup\{+\infty\}$).\vspace{.8pc}

\noindent We say that the sequence $\{\varphi\}_{n\geq 1}$
converges in the Mosco sense to $\varphi$, denoted by
$\varphi_n\stackrel{M}{\longrightarrow}\varphi$ if and only if the
following two conditions hold:
\begin{enumerate}
\renewcommand\labelenumi{(\arabic{enumi})}
\leftskip .1pc
\item for every $y\in Y$ and every sequence $\{y_n\}_{n\geq
1}\subseteq Y$ such that $y_n\stackrel{w}{\longrightarrow}y$ in
$Y$, we have that
\begin{equation*}
\hskip -1.25pc \varphi(y) \leq
\liminf_{n\rightarrow+\infty}\varphi_n(y_n);
\end{equation*}
\item for every $y\in Y$, there exists a sequence $\{y_n\}_{n\geq
1}\subseteq Y$ such that $y_n\longrightarrow y$ in $Y$ and
$\varphi_n(y_n)\longrightarrow\varphi(y)$.
\end{enumerate}

Further analysis of these two notions can be found in \cite{HPa}.

\section{The nonsmooth reduction method} \label{reduct}

In this section we extend the Castro, Lazer and Thews reduction
method to the present nonsmooth setting. By $2^{*}$ we denote the
Sobolev critical exponent defined by
\begin{equation*}
2^{*} \stackrel{d f}{=} \begin{cases}
\displaystyle \frac{2N}{N-2}, &\hbox{if} \quad N>2 \\[.4pc]
+\infty, &\hbox{if} \quad N\leq 2. \end{cases}
\end{equation*}
Our hypotheses on the nonsmooth potential $j$ are the
following:\vspace{.5pc}

\noindent ${\rm H(j)}\hbox{\rm :}\ j\hbox{\rm :}\
Z\times\mathbb{R}\longrightarrow \mathbb{R}$ is a function, such
that

\begin{enumerate}
\renewcommand\labelenumi{(\roman{enumi})}
\leftskip .6pc
\item for all $\zeta\in\mathbb{R}$, the function $Z\ni z\longmapsto
j(z,\zeta)\in\mathbb{R}$ is measurable and for almost all $z\in
Z$, we have that $j(z,0)=0$;

\item for almost all $z\in Z$, the function
$\mathbb{R}\ni\zeta\longmapsto j(z,\zeta)\in\mathbb{R}$ is locally
Lipschitz;

\item for almost all $z\in Z$, all $\zeta\in\mathbb{R}$ and all
$u\in\partial j(z,\zeta)$, we have that $|u|\leq
a_1(z)+c_1|\zeta|^{r-1}$, with $a_1\in L^{\infty}(Z)$, $c_1>0$ and
$1\leq r< 2^*$;

\item $\lim_{|\zeta|\rightarrow+\infty}[
u(\zeta)\zeta-2j(z,\zeta)]=-\infty$ uniformly for almost all $z\in
Z$ and all $u(\zeta)\in\partial j(z,\zeta)$;

\item there exists $l\in L^{\infty}(Z)$, such that for almost all $z\in
Z$, we have $\l(z)\leq\lambda_{k+1}-\lambda_k$ with strict
inequality on a set of positive measure and for almost all $z\in
Z$, all $\zeta_1,\zeta_2\in\mathbb{R}$ with $\zeta_1\not=\zeta_2$
and all $v_1\in\partial j(z,\zeta_1)$, $v_2\in\partial
j(z,\zeta_2)$, we have $\frac{v_1-v_2}{\zeta_1-\zeta_2}\leq l(z)$;

\item there exist $\beta\in L^{\infty}(Z)_-$ and $\delta_0>0$ such that
for some integer $m\in[1,k]$ and for almost all $z\in Z$, we have
$\beta(z)\leq\lambda_m-\lambda_k$, with strict inequality on a set
of positive measure and for almost all $z\in Z$ and all
$\zeta\in\mathbb{R}$, such that $|\zeta|\leq\delta_0$, we have
$\lambda_{m-1}-\lambda_k\leq
\frac{2j(z,\zeta)}{\zeta^2}\leq\beta(z)$;

\item $0\leq\liminf_{|\zeta|\rightarrow+\infty}
\frac{2j(z,\zeta)}{\zeta^2}
\leq\limsup_{|\zeta|\rightarrow+\infty}
\frac{2j(z,\zeta)}{\zeta^2} \leq\gamma(z)$ uniformly for almost
all $z\in Z$, with $\gamma\in L^{\infty}(Z)_+$ and
$\gamma(z)\leq\lambda_{k+1}-\lambda_k$ for almost all $z\in Z$,
with strict inequality on a set of positive measure.
\end{enumerate}

\setcounter{defin}{0}
\begin{remark}
{\rm Hypothesis {\rm H(j)(vi)} implies that we have double
resonance at the origin. The resonance is complete from below and
incomplete from above. The same double resonance situation at
infinity is implied by hypothesis {\rm H(j)(vii)}. So hypotheses
{\rm H(j)(vi)} and {\rm H(j)(vii)} together provide the
double-double resonance character of our problem. Also note that
hypothesis {\rm H(j)(v)} permits only downward discontinuities of
the derivative of the potential function $j(z,\cdot)$. Recall that
for almost all $z\in Z$, the derivative of $j(z,\cdot)$ exists
almost everywhere (Rademacher theorem).}
\end{remark}

Let $\varphi\hbox{\rm :}\ H^1_0(Z)\longrightarrow\mathbb{R}$ be
the energy functional defined by
\begin{equation*}
\varphi(x) \stackrel{d f}{=} \frac{1}{2}\|\nabla x\|_2^2
-\frac{\lambda_k}{2}\|x\|_2^2 -\int_Zj(z,x(z))\hbox{d} z
\quad\forall x\in H^1_0(Z)
\end{equation*}
(by $\|\cdot\|_p$ we denote the norm of $L^p(Z)$). We know that
$\varphi$ is locally Lipschitz (see p.~313 of \cite{HPb}). Since
$k\geq 1$ is fixed, for what follows we set
\begin{equation*}
\bar{H} \stackrel{d f}{=} \bar{H}_k = \bigoplus_{i=1}^{k-1}
E(\lambda_i) \quad\hbox{and}\quad \hat{H} \stackrel{d f}{=}
\hat{H}_k = \overline{\bigoplus_{i=k+1}^{\infty} E(\lambda_i)}.
\end{equation*}
We have that
\begin{equation*}
H^1_0(Z) = \bar{H}\oplus E(\lambda_k)\oplus\hat{H}.
\end{equation*}
Also we set
\begin{equation*}
\bar{H}_0 \stackrel{d f}{=} \bar{H}\oplus E(\lambda_k)
 = \bigoplus_{i=1}^k E(\lambda_i)
\end{equation*}
and for $u\in\bar{H}_0$, we consider the following minimization
problem: \setcounter{equation}{0}
\begin{equation} \label{eq3}
\inf_{v\in\hat{H}}\varphi(u+v).
\end{equation}

Since we do not identify $H^1_0(Z)$ with its dual, we have that
\begin{equation*}
(H^1_0(Z))^* = H^{-1}(Z) = \bar{H}_0^*\oplus\hat{H}^*.
\end{equation*}

We start with a simple lemma which is needed in what follows.

\begin{lemma} \label{lem2}
If $n\geq 1$ and $\beta\in L^{\infty}(Z)_+${\rm ,} with
\begin{equation*}
\beta(z)\leq\lambda_{n+1} \ \ \hbox{for a.a.}\ \ z\in Z
\end{equation*}
and the inequality is strict on a set of positive measure{\rm ,}
then there exists $\xi_1>0${\rm ,} such that
\begin{equation*}
\|\nabla x\|_2^2-\int_Z\beta(z)|x(z)|^2\hbox{\rm d} z \geq
\xi_1\|\nabla x\|_2^2 \quad\forall x\in \hat{H}_n.
\end{equation*}
\end{lemma}

\begin{proof}
Let
\begin{equation*}
\theta(x) \stackrel{d f}{=} \|\nabla x\|_2^2
-\int_Z\beta(z)|x(z)|^2 \hbox{d} z \quad\forall x\in \hat{H}_n.
\end{equation*}
By virtue of the variational characterization of the eigenvalues,
we have that $\theta\geq 0$. Suppose that the lemma is not true.
Because of the positive homogeneity of $\theta$, we can find a
sequence $\{x_m\}_{m\geq 1}\subseteq \hat{H}_n$, such that
$\|\nabla x_m\|_2=1$ for $m\geq 1$ and $\theta(x_m)\searrow 0$. By
passing to a subsequence if necessary, we may assume that
\begin{align*}
x_m\stackrel{w}{\longrightarrow} x_0 &\quad \hbox{in}\ H^1_0(Z),\\[.4pc]
x_m\longrightarrow x_0               &\quad \hbox{in}\ L^2(Z),
\end{align*}
with some $x_0\in \hat{H}_n$. Since the norm in a Banach space is
weakly lower semicontinuous, in the limit as
$m\rightarrow+\infty$, we obtain
\begin{equation*}
\theta(x_0) = \|\nabla x_0\|_2^2-\int_Z \beta(z)|x_0(z)|^2
\hbox{d} z \leq 0,
\end{equation*}
so
\begin{equation*}
\|\nabla x_0\|_2^2 \leq \int_Z\beta(z)|x_0(z)|^2 \hbox{d} z  \leq
\lambda_{n+1}\|x_0\|_2^2
\end{equation*}
and since $x_0\in\hat{H}_n$, from the variational characterization
of $\lambda_{n+1}$, we have that
\begin{equation} \label{eq4}
\|\nabla x_0\|_2^2 = \lambda_{n+1}\| x_0\|_2^2.
\end{equation}

If $x_0=0$, then taking into account that $\theta(x_m)\to 0$ we
would have that $\|\nabla x_m\|_2\to 0$. Because $\|\nabla
x_m\|_2=1$ for $m\geq 1$, this is not possible, so $x_0\not= 0$.
From~(\ref{eq4}), it follows that $x_0\in E(\lambda_{n+1})$. Then,
from the hypothesis that $\beta(z)<\lambda_{n+1}$ on a set of
positive measure and from the unique continuity property of
$E(\lambda_{n+1})$, we have that
\begin{equation*}
\|\nabla x_0\|_2^2 \leq \int_Z\beta(z)|x_0(z)|^2\hbox{d} z <
\lambda_{n+1}\|x_0\|_2^2,
\end{equation*}
a contradiction to~(\ref{eq4}).\hfill $\Box$
\end{proof}

The next proposition essentially extends the Castro--Lazer--Thews
reduction method to a nonsmooth setting.

\setcounter{defin}{2}
\begin{proposition}\label{prop3}$\left.\right.$\vspace{.5pc}

\noindent If hypotheses ${\rm H(j)}$ hold{\rm ,} then there exists
a continuous map $\vartheta\hbox{\rm :}\
\bar{H}_0\longrightarrow\hat{H}${\rm ,} such that for every
$u\in\bar{H}_0${\rm ,} we have
\begin{equation*}
\inf_{v\in\hat{H}}\varphi(u+v) = \varphi(u+\vartheta(u))
\end{equation*}
and $\vartheta(u)\in\hat{H}$ is the unique solution of the
operator inclusion
\begin{equation*}
0 \in p_{\hat{H}^*}\partial\varphi(u+v),
\end{equation*}
with $u\in\bar{H}_0$ fixed{\rm ,} where $p_{\hat{H}^*}$ is the
orthogonal projection on $\hat{H}^*=[\bar{H}^*_0]^{\perp}$.
\end{proposition}

\begin{proof}
For a fixed $u\in\bar{H}_0$, let $\varphi_u\hbox{\rm :}\
H^1_0(Z)\longrightarrow\mathbb{R}$ be defined by
\begin{equation*}
\varphi_u(w) \stackrel{d f}{=} \varphi(u+w) \quad\forall w\in
H^1_0(Z).
\end{equation*}
For every $w,h\in H^1_0(Z)$, we have that
\begin{align*}
\varphi_u^0(w;h) &= \limsup_{\substack{
w'\,\rightarrow \,w\\[.4pc]
t\,\searrow \,0}}
\frac{\varphi_u(w'+th)-\varphi_u(w')}{t}\\[.4pc]
&=\limsup_{\substack{
w'\,\rightarrow \,w\\[.4pc]
t\,\searrow \,0 }} \frac{\varphi(u+w'+th)-\varphi(u+w')}{t} =
\varphi^0(u+w;h),
\end{align*}
so
\begin{equation} \label{eq5}
\partial\varphi_u(w) = \partial\varphi(u+w) \quad\forall w\in H^1_0(Z).
\end{equation}

Let $\hat{i}\hbox{\rm :}\ \hat{H}\longrightarrow H^1_0(Z)$ be the
inclusion map and let $\hat{\varphi}_u\hbox{\rm :}\
\hat{H}\longrightarrow\mathbb{R}$ be defined by
\begin{equation*}
\hat{\varphi}_u(v)=\varphi(u+v) \quad\forall v\in\hat{H}.
\end{equation*}
We have that $\varphi_u\circ\hat{i}=\hat{\varphi}_u$ and so
\begin{equation} \label{eq6}
\partial(\varphi_u\circ\hat{i})(v) =
\partial\hat{\varphi}_u(v) \quad\forall v\in\hat{H}.
\end{equation}

But from the chain rule of Clarke (p.~45--46 of \cite{Cl}), we
have
\begin{equation*}
\partial(\varphi_u\circ\hat{i})(v) \subseteq
p_{\hat{H}^*}\partial\varphi_u(\hat{i}(v)) \quad\forall
v\in\hat{H},
\end{equation*}
since $\hat{i}^*=p_{\hat{H}^*}$. Hence, from~(\ref{eq5})
and~(\ref{eq6}), it follows that
\begin{equation} \label{eq7}
\partial\hat{\varphi}_u(v) \subseteq
p_{\hat{H}^*}\partial\varphi_u(\hat{i}(v)) =
p_{\hat{H}^*}\partial\varphi(u+v) \quad\forall v\in\hat{H}.
\end{equation}

Now we have
\begin{equation} \label{eql_60}
x^* = A(x)-\lambda_k x-h \quad\forall x\in H^1_0(Z),
x^*\in\partial\varphi(x),
\end{equation}
with $A\in\mathcal{L}( H^1_0(Z),H^{-1}(Z))$ being defined by
\begin{equation*}
\langle A(x),y\rangle_{H^1_0(Z)} \stackrel{d f}{=} \int_Z(\nabla
x(z),\nabla y(z))_{\mathbb{R}^N}\hbox{d} z \quad\forall x,y\in
H^1_0(Z)
\end{equation*}
and $h\in L^{r'}(Z)$ (where $\frac{1}{r}+\frac{1}{r'}=1$), such
that $h(z)\in\partial j(z,x(z))$ for almost all $z\in Z$ (see
p.~83 of \cite{Cl}). So for any $v_1,v_2\in\hat{H}$ and
$x^*_1\in\partial\hat{\varphi}_u(v_1)$,
$x^*_2\in\partial\hat{\varphi}_u(v_2)$, we have
\begin{equation*}
x_i^* = p_{\hat{H}^*}A(u+v_i)-\lambda_k v_i-p_{\hat{H}^*}h_i
\quad\hbox{for}\ i\in\{1,2\},
\end{equation*}
where $h_i\in L^{r'}(Z)\subseteq H^{-1}(Z)$ (recall $r<2^*$) is
such that $h_i(z)\in\partial j(z,(u+v_i)(z))$ for almost all $z\in
Z$ and $i\in\{1,2\}$. Since $p^*_{\hat{H}^*}=\hat{i}$, we have
that
\begin{align*}
&\langle p_{\hat{H}^*}(A(u+v_1)-A(u+v_2)),
v_1-v_2\rangle_{\hat{H}}\\[.4pc]
&\quad\, = \langle A(u+v_1)-A(u+v_2),
\hat{i}(v_1)-\hat{i}(v_2)\rangle_{H^1_0(Z)}\\[.4pc]
&\quad\, = \langle A(u+v_1)-A(u+v_2), v_1-v_2\rangle_{H^1_0(Z)} =
\|\nabla v_1-\nabla v_2\|_2^2.
\end{align*}

By hypothesis H(j)(v), we have that
\begin{equation} \label{eq8}
\frac{h_1(z)-h_2(z)}{v_1(z)-v_2(z)} \leq l(z) \quad\hbox{a.e. on}\
\{v_1\not= v_2\}.
\end{equation}
So, from~(\ref{eql_60}) and~(\ref{eq8}), it follows that
\begin{align*}
&\langle x_1^*-x_2^*, v_1-v_2\rangle_{\hat{H}}\\[.4pc]
&\quad =\|\nabla v_1-\!\nabla v_2\|_2^2 -\lambda_k\| v_1-v_2\|_2^2
-\int_Z(h_1(z)-h_2(z))(v_1(z)-v_2(z))\hbox{d} z\\[.4pc]
&\quad \geq \|\nabla v_1-\!\nabla v_2\|_2^2 -\lambda_k\|
v_1-v_2\|_2^2 -\int_Zl(z)|v_1(z)-v_2(z)|^2\hbox{d} z.
\end{align*}

By hypothesis H(j)(v), we know that
\begin{equation*}
l(z) \leq \lambda_{k+1}-\lambda_k
\quad\hbox{for a.a.}\ z\in Z,
\end{equation*}
with strict inequality on a set of positive  measure. So we can
apply Lemma~\ref{lem2} (with $\beta(z)=l(z)+\lambda_k$) and obtain
$\xi_1>0$, such that
\begin{equation*}
\langle x_1^*-x_2^*,v_1-v_2\rangle_{\hat{H}} \geq \xi_1\|\nabla
v_1-\nabla v_2\|_2^2.
\end{equation*}
So the multifunction $v\longmapsto\partial\hat{\varphi}_u(v)$ is
strongly monotone in the dual pair $(\hat{H},\hat{H}^*)$. Hence
the function $\hat{H}\ni
v\longmapsto\hat{\varphi}_u(v)\in\mathbb{R}$ is strongly convex,
i.e. the function $\hat{H}\ni v\longmapsto\hat{\varphi}_u(v)
-\frac{\xi_1}{2}\|v\|^2_{H^1_0(Z)}\in\mathbb{R}$ is convex (see
p.~37 of \cite{Cl}).

Let $v\in\hat{H}$, $x^*\in\partial\hat{\varphi}_u(v)$ and
$y^*\in\partial\hat{\varphi}_u(0)$. From the previous
considerations, we have
\begin{align*}
\langle x^*,v\rangle_{\hat{H}} &= \langle x^*-y^*,
v\rangle_{\hat{H}} +\langle y^*,v\rangle_{\hat{H}}\\[.4pc]
&\geq\xi_1\|\nabla v\|_2^2-\xi_2\|y^*\|_{H^{-1}(Z)} \|\nabla
v\|_2,
\end{align*}
for some $\xi_2>0$, so the multifunction
$v\longmapsto\partial\varphi_u(v)$ is coercive.

The multifunction $v\longmapsto\partial\hat{\varphi}_u(v)$ is
maximal monotone (since $\hat{\varphi}_u$ is convex). But a
maximal monotone, coercive operator is surjective (see p.~322 of
\cite{HPa}). Thus, we can find $v_0\in\hat{H}$, such that
\begin{equation*}
0 \in \partial\hat{\varphi}_u(v_0) \quad\hbox{and}\quad
\inf_{v\in\hat{H}}\varphi(u+v) = \varphi(u+v_0).
\end{equation*}
Because of the strong convexity of $\hat{\varphi}_u$, we infer
that the minimizer $v_0\in\hat{H}$ is unique.

Therefore we can define a map $\vartheta\hbox{\rm :}\
\bar{H}_0\longrightarrow\hat{H}$ which to each fixed
$u\in\bar{H}_0$ assigns the unique solution $v_0\in\hat{H}$ of the
minimization problem~(\ref{eq3}). Then from~(\ref{eq7}), we have
\begin{equation*}
0 \in \partial\hat{\varphi}_u(\vartheta(u)) \subseteq
p_{\hat{H}^*}\partial\varphi(u+\vartheta(u))
\end{equation*}
and
\begin{equation*}
\inf_{v\in\hat{H}}\varphi(u+v) = \varphi(u+\vartheta(u)).
\end{equation*}

Finally, we have to show that $\vartheta$ is continuous. To this
end suppose that $u_n\rightarrow u$ in $\bar{H}_0$. If
$v_n\longrightarrow v$ in $\hat{H}$, we have
$\hat{\varphi}_{u_n}(v_n)\longrightarrow\hat{\varphi}_u(v)$ (in
fact it is easy to see that
$\hat{\varphi}_{u_n}\longrightarrow\hat{\varphi}_u$ in
$C(\hat{H})$). On the other hand, if
$v_n\stackrel{w}{\longrightarrow}v$ in $\hat{H}$, by virtue of the
weak lower semicontinuity of the norm in a Banach space and the
compactness of the embedding $H^1_0(Z)\subseteq L^2(Z)$, we have
\begin{equation*}
\hat{\varphi}_u(v) \leq \liminf_{n\rightarrow+\infty}
\hat{\varphi}_{u_n}(v_n)
\end{equation*}
(recall the definition of $\varphi$). If follows that
$\hat{\varphi}_{u_n}\stackrel{M}{\longrightarrow} \hat{\varphi}_u$
and so by virtue of Theorem~5.6.9 of p.~766 of  \cite{HPa}, we
have that
\begin{equation*}
\hbox{Gr}\,\partial\hat{\varphi}_{u_n}
\stackrel{K_{s,s}}{\longrightarrow}
\hbox{Gr}\,\partial\hat{\varphi}_u \quad\hbox{as}\
n\rightarrow+\infty.
\end{equation*}
Because $0\in\partial\hat{\varphi}_u(\vartheta(u))$, we can find
$v_n^*\in\partial\hat{\varphi}_{v_n}(v_n)$, for $n\geq 1$, such
that $v_n\rightarrow \vartheta(u)$ in $\hat{H}$ and
$v_n^*\longrightarrow 0$ in $\hat{H}^*$. Recall that
$0\in\partial\hat{\varphi}_{u_n}(\vartheta(u_n))$, from the strong
monotonicity of $\partial\hat{\varphi}_{u_n}$, we have that
\begin{equation*}
\langle v_n^*,v_n-\vartheta(u_n)\rangle_{\hat{H}} \geq
\xi_1\|v_n-\vartheta(u_n)\|_{H^1_0(Z)}^2,
\end{equation*}
so
\begin{align*}
&\|\vartheta(u)-\vartheta(u_n)\|_{H^1_0(Z)}\\[.4pc]
&\quad\,\leq
\frac{1}{\xi_1}\|v_n^*\|_{H^{-1}(Z)}+\|v_n-\vartheta(u)\|_{H_0^1(Z)}
\longrightarrow 0 \quad\hbox{as}\ n\rightarrow+\infty
\end{align*}
and thus we have proved that $\vartheta$ is continuous.\hfill
$\Box$
\end{proof}

Using Proposition~\ref{prop3}, we can define the map
$\bar{\varphi}\hbox{\rm :}\ \bar{H}_0\longrightarrow\mathbb{R}$,
by
\begin{equation} \label{eql_70}
\bar{\varphi}(u) = \varphi(u+\vartheta(u)) \quad\forall
u\in\bar{H}_0.
\end{equation}
Note that from the definition of $\vartheta$, for all
$u,h\in\bar{H}_0$, we have that
\begin{align*}
\bar{\varphi}(u+h)-\bar{\varphi}(u) &=
\varphi(u+h+\vartheta(u+h))-\varphi(u+\vartheta(u))\\[.4pc]
&\leq \varphi(u+h+\vartheta(u))-\varphi(u+\vartheta(u)).
\end{align*}
Similarly from the definition of $\vartheta$, for all
$u,h\in\bar{H}_0$, we obtain
\begin{align*}
\bar{\varphi}(u)-\bar{\varphi}(u+h) &=
\varphi(u+\vartheta(u))-\varphi(u+h+\vartheta(u+h))\\[.4pc]
&\leq \varphi(u+\vartheta(u+h))-\varphi(u+h+\vartheta(u+h)).
\end{align*}
If follows that $\bar{\varphi}$ is locally Lipschitz (since
$\varphi$ is).

Now we will show that
\begin{equation} \label{eq10}
\partial\bar{\varphi}(u) \subseteq
p_{\bar{H}_0^*}\partial\varphi(u+\vartheta(u)) \quad\forall
u\in\bar{H}_0.
\end{equation}
First for all $u,h\in\bar{H}_0$, we have
\begin{align*}
\bar{\varphi}^0(u;h) &= \limsup_{\substack{
u'\,\rightarrow \,u\\[.4pc]
t\,\searrow \,0 }}
\frac{\bar{\varphi}(u'+th)-\bar{\varphi}(u')}{t}\\[.4pc]
&= \limsup_{\substack{
u'\,\rightarrow \,u\\[.4pc]
t\,\searrow \,0}} \frac{\varphi(u'+th+\vartheta(u'+th))
-\varphi(u'+\vartheta(u'))}{t}\\[.4pc]
&\leq \limsup_{\substack{
u'\,\rightarrow \,u\\[.4pc]
t\,\searrow \,0}} \frac{\varphi(u'+th+\vartheta(u'))
-\varphi(u'+\vartheta(u'))}{t}\\[.4pc]
&\leq \varphi^0 (u+\vartheta(u);h).
\end{align*}
Denoting by $\bar{i}_0\hbox{\rm :}\ \bar{H}_0\longrightarrow
H^1_0(Z)$ the inclusion map and noting that
$\bar{i}_0^*=p_{\bar{H}_0^*}$, for all $u,h\in\bar{H}_0$, we have
that
\begin{align*}
&\bar{\varphi}^0(u;h) \leq \varphi^0
(u+\vartheta(u);\bar{i}_0(h))\\[.4pc]
&\quad\,= \sup_{u^*\in\partial\varphi(u+\vartheta(u))} \langle
u^*,\bar{i}_0(h) \rangle_{H^1_0(Z)} =
\sup_{u^*\in\partial\varphi(u+\vartheta(u))} \langle
p_{\bar{H}_0^*}(u^*),h \rangle_{\bar{H}_0}.
\end{align*}
Suppose that $u_0^*\in\partial\bar{\varphi}(u)$. From the
definition of the Clarke directional derivative, we have
\begin{equation*}
\langle u_0^*,h\rangle_{\bar{H}_0} \leq \bar{\varphi}^0(u;h)
\quad\forall h\in\bar{H}_0,
\end{equation*}
so
\begin{equation*}
\langle u_0^*,h\rangle_{\bar{H}_0} \leq
\sup_{u^*\in\partial\varphi(u+\vartheta(u))} \langle
p_{\bar{H}_0^*}(u^*),h \rangle_{\bar{H}_0} \qquad\forall
h\in\bar{H}_0
\end{equation*}
and thus
\begin{equation*}
u_0^* \in p_{\bar{H}_0^*} \partial\varphi(u+\vartheta(u)).
\end{equation*}
Therefore, we obtain~(\ref{eq10}).

Next let $\psi=-\bar{\varphi}$. Then $\psi$ is locally Lipschitz
on the finite dimensional space $\bar{H}_0$. In the next section
working with $\psi$ and using Theorem~\ref{thm1}, we prove a
multiplicity theorem for problem~(HVI).

\section{Existence of multiple solutions} \label{exist}

As $m\leq k$ are fixed (see hypothesis H(j)(vi)), let us put
\begin{equation*}
Y \stackrel{d f}{=} \bigoplus_{i=1}^{m-1}E(\lambda_i)
\quad\hbox{and}\quad V \stackrel{d f}{=} \bigoplus_{i=m}^k
E(\lambda_i).
\end{equation*}
We have
\begin{equation*}
\bar{H}_0 = \bar{H}\oplus E(\lambda_k) = Y\oplus V.
\end{equation*}
The next proposition shows that $\psi=-\bar{\varphi}$ satisfies
the local linking condition (see Theorem~\ref{thm1}).

\setcounter{defin}{0}
\begin{proposition} \label{prop4}$\left.\right.$\vspace{.5pc}

\noindent If hypotheses {\rm H(j)} hold{\rm ,} then there exists
$\delta>0$ such that
\begin{equation*}
\begin{cases}
\psi(u)\leq 0 &\hbox{if} \ u\in V,\|u\|_{H^1_0(Z)}\leq\delta\\[.4pc]
\psi(u)\geq 0 &\hbox{if} \ u\in Y,\|u\|_{H^1_0(Z)}\leq\delta
\end{cases}.
\end{equation*}
\end{proposition}

\begin{proof}
Because $Y$ is finite dimensional, all norms are equivalent and
so we can find $M_1>0$, such that
\begin{equation*}
\sup_{z\in Z} |u(z)| \leq M_1\|u\|_{H^1_0(Z)} \quad\forall u\in Y.
\end{equation*}
Let $\beta\in L^{\infty}(Z)$ and $\delta_0>0$ as in hypothesis
H(j)(vi). Thus, if $\delta'\stackrel{df}{=}\frac{\delta_0}{M_1}$,
we have
\begin{equation*}
\sup_{z\in Z}|u(z)| \leq \delta_0 \quad\forall u\in Y,
\|u\|_{H^1_0(Z)}\leq\delta'.
\end{equation*}
By virtue of hypothesis H(j)(vi) and the definition of
$\vartheta$, for all $u\in Y$ with $\|u\|_{H^1_0(Z)}\leq\delta'$,
we have
\begin{align*}
\psi(u) &= -\bar{\varphi}(u) \geq -\frac{1}{2}\|\nabla u\|_2^2
+\frac{\lambda_k}{2}\|u\|_2^2
+\int_Zj(z,u(z))\hbox{d} z\\[.4pc]
&\geq -\frac{1}{2}\|\nabla u\|_2^2 +\frac{\lambda_k}{2}\|u\|_2^2
+\frac{\lambda_{m-1}-\lambda_k}{2}\|u\|_2^2\\[.4pc]
&= -\frac{1}{2}\|\nabla u\|_2^2 +\frac{\lambda_{m-1}}{2}\|u\|_2^2
\geq 0.
\end{align*}

Also for all $u\in V$, we have \setcounter{equation}{0}
\begin{align} \label{eq11}
\psi(u) &= -\varphi(u+\vartheta(u))\\[.4pc]
&= -\frac{1}{2}\|\nabla(u+\vartheta(u))\|_2^2
+\frac{\lambda_k}{2}\|u+\vartheta(u)\|_2^2
+\!\int_Zj(z,(u+\vartheta(u))(z))\hbox{d} z.\nonumber
\end{align}

From hypothesis H(j)(vi), we have that
\begin{equation*}
j(z,\zeta) \leq \frac{1}{2}\beta(z)\zeta^2 \quad \hbox{for a.a.}\
z\in Z \ \hbox{and all}\ |\zeta|\leq\delta_0,
\end{equation*}
while by virtue of hypothesis H(j)(iii) and the Lebourg mean value
theorem (see p.~41 of \cite{Cl} or \cite{L}), we have that
\begin{equation*}
j(z,\zeta) \leq c_2|\zeta|^{\eta} \quad \hbox{for a.a.}\ z\in Z, \
\hbox{all}\ |\zeta|>\delta_0,
\end{equation*}
with some $2<\eta\leq 2^*$ and $c_2>0$. So finally we can say that
\begin{equation*}
j(z,\zeta) \leq \frac{1}{2}\beta(z)\zeta^2+c_3|\zeta|^{\eta}
\quad\hbox{for a.a.}\ z\in Z,\ \hbox{all} \ \zeta\in\mathbb{R}
\end{equation*}
with $c_3\stackrel{d f}{=} c_2+\frac{1}{2}\|\beta\|_{\infty}>0$.
Using this in~(\ref{eq11}), we obtain
\begin{align} \label{eq12}
\psi(u) &\leq -\frac{1}{2}\|\nabla(u+\vartheta(u))\|_2^2
+\frac{\lambda_k}{2}\|u+\vartheta(u)\|_2^2\\[.4pc]
&\quad\, +\frac{1}{2}\int_Z\beta(z)|(u+\vartheta(u))(z)|^2\hbox{d}
z +c_3\|u+\vartheta(u)\|_{\eta}^{\eta} \quad\forall u\in
V.\nonumber
\end{align}

Note that
\begin{equation*}
\lambda_{m-1} \leq \beta(z)+\lambda_k \leq \lambda_{m}
\quad\hbox{for almost all}\ z\in Z,
\end{equation*}
with the second inequality strict on a set of positive measure.
Since $u\in V$ and $\vartheta(u)\in\hat{H}$,
$u+\vartheta(u)\in\hat{H}_{m-1}$ and we can apply Lemma~\ref{lem2}
with $\beta+\lambda_k$ to obtain $\xi_1>0$ such that
\begin{align*}
&\|\nabla(u+\vartheta(u))\|_2^2
-\int_Z(\beta(z)+\lambda_k)|(u+\vartheta(u))(z)|^2 \hbox{d} z\\[.4pc]
&\quad\,\geq \xi_1\|\nabla(u+\vartheta(u))\|_2^2 \quad\forall v\in
V.
\end{align*}
Using this inequality in~(\ref{eq12}), we have
\begin{equation*}
\psi(u) \leq -\frac{1}{2}\xi_1\|\nabla(u+\vartheta(u))\|_2^2
+c_4\|\nabla(u+\vartheta(u))\|_2^{\eta} \quad\forall v\in V,
\end{equation*}
for some $c_4>0$. Here we have used the Sobolev embedding theorem
since $\eta\leq 2^*$ and the Poincar\'e inequality. Because
$2<\eta$, $\vartheta(0)=0,$ we can find $\delta''>0$, such that
\begin{equation*}
\psi(u)\leq 0 \quad\forall u\in V, \|u\|_{H^1_0(Z)}\leq\delta''.
\end{equation*}
Finally let $\delta\stackrel{d f}{=}\min\{\delta',\delta''\}$ to
finish the proof of the proposition.\hfill $\Box$
\end{proof}

Since we aim in applying Theorem~\ref{thm1}, we need to show that
$\psi$ satisfies the nonsmooth C-condition. To establish this for
$\psi$, we show first that $\varphi$ satisfies the nonsmooth\break
C-condition.

\begin{proposition} \label{prop5}$\left.\right.$\vspace{.5pc}

\noindent If hypotheses {\rm H(j)} hold{\rm ,} then $\varphi$
satisfies the nonsmooth ${\rm C}$-condition.
\end{proposition}

\begin{proof}
Let $\{x_n\}_{n\geq 1}\subseteq H^1_0(Z)$ be a sequence, such that
\begin{equation*}
\varphi(x_n)\longrightarrow c
\quad\hbox{and}\quad
(1+\|x_n\|_{H^1_0(Z)})m_{\varphi}(x_n)\longrightarrow 0
\ \hbox{as}\ n\rightarrow+\infty.
\end{equation*}
Since the norm functional is weakly lower semicontinuous and
$\partial\varphi(x_n)\subseteq H^{-1}(Z)$ is weakly compact, from
the Weierstrass theorem, we know that there exists
$x_n^*\in\partial\varphi(x_n)$, such that
$\|x_n^*\|_{H^{-1}(Z)}=m_{\varphi}(x_n)$ for all $n\geq 1$. We
know that
\begin{equation*}
x_n^* =  A(x_n)-\lambda_kx_n-h_n \quad\forall n\geq 1,
\end{equation*}
with $A\in\mathcal{L}( H^1_0(Z),H^{-1}(Z))$ being the maximal
monotone operator defined by
\begin{equation*}
\langle A(x),y\rangle_{H^1_0(Z)} \stackrel{d f}{=} \int_Z(\nabla
x(z),\nabla y(z))_{\mathbb{R}^N}\hbox{d} z \qquad\forall x,y\in
H^1_0(Z)
\end{equation*}
and $h_n\in L^{r'}\!(Z)$ (with $\frac{1}{r}+\frac{1}{r'}=1$) is
such that $h_n(z)\in\partial j(z,x_n(z))$ for almost all $z\in Z$
(see p.~80 of \cite{Cl}). From the choice of the sequence
$\{x_n\}_{n\geq 1}\subseteq H^1_0(Z)$, we have
\begin{align}
&| \langle x_n^*,x_n\rangle_{H^1_0(Z)} - 2\varphi(x_n)+2c|\nonumber\\[.4pc]
&\quad\, \leq \|x_n\|_{H^1_0(Z)}\|x_n^*\|_{H^{-1}(Z)}+2|\varphi(x_n)-c|\nonumber\\[.4pc]
\label{eq13} &\quad\, \leq (1+\|x_n\|_{H^1_0(Z)})m_{\varphi}(x_n)
+2|\varphi(x_n)-c| \longrightarrow 0 \quad \hbox{as}\
n\rightarrow+\infty.
\end{align}
Note that
\begin{equation*}
\langle x_n^*,x_n\rangle_{H^1_0(Z)} = \langle
A(x_n),x_n\rangle_{H^1_0(Z)}-\lambda_k\|x_n\|_2^2 -\int_Z h_n(z)
x_n(z) \hbox{d} z.
\end{equation*}
Then using~(\ref{eq13}), we obtain
\begin{align} \label{eq14}
&\int_Z(h_n(z)x_n(z)-2j(z,x_n(z)))\hbox{d} z\nonumber\\[.4pc]
&\quad\, = 2\varphi(x_n)-\langle x_n^*,x_n\rangle_{H^1_0(Z)}
\longrightarrow 2c \quad \hbox{as}\ n\rightarrow+\infty.
\end{align}
We claim that $\{x_n\}_{n\geq 1}\subseteq H^1_0(Z)$ is bounded.
Suppose that this is not the case. By passing to a subsequence if
necessary, we may assume that $\|x_n\|_{H^1_0(Z)}\longrightarrow
+\infty$. Let $y_n\stackrel{d f}{=}\frac{x_n}{\|x_n\|_{H^1_0(Z)}}$
for $n\geq 1$. Since $\|y_n\|_{H^1_0(Z)}=1$, passing to a
subsequence if necessary, we may assume that
\begin{align*}
&y_n\stackrel{w}{\longrightarrow}y \quad \hbox{in}\ H^1_0(Z),\\[.4pc]
&y_n\longrightarrow y              \quad \hbox{in}\ L^2(Z),\\[.4pc]
&y_n(z)\longrightarrow y(z)        \quad \hbox{for a.a.}\ z\in Z
\end{align*}
and
\begin{align*}
&|y_n(z)|\leq k(z) \quad \hbox{for a.a.}\ z\in Z\ \hbox{and all}\
n\geq 1,
\end{align*}
with some $k\in L^2(Z)$. Because of hypothesis H(j)(vii), for a
given $\varepsilon>0$, we can find $M_2=M_2(\varepsilon)>0$ such
that for almost all $z\in Z$ and all $|\zeta|\geq M_2$, we have
\begin{equation*}
j(z,\zeta) \leq \frac{1}{2}(\gamma(z)+\varepsilon)\zeta^2.
\end{equation*}

Also from hypothesis H(j)(iii) and the Lebourg mean value theorem,
for almost all $z\in Z$ and all $|\zeta|<M_2$, we have
$|j(z,\zeta)|\leq\xi_2$ for some $\xi_2>0$. So we can say that for
almost all $z\in Z$ and all $\zeta\in\mathbb{R}$, we have
\begin{equation*}
j(z,\zeta) \leq \frac{1}{2}(\gamma(z)+\varepsilon)\zeta^2+\xi_2.
\end{equation*}
Then for every $n\geq 1$, we have
\begin{align*}
\frac{\varphi(x_n)}{\|x_n\|_{H^1_0(Z)}^2} &=\frac{1}{2}\|\nabla
y_n\|_2^2 -\frac{\lambda_k}{2}\| y_n\|_2^2
-\int_Z\frac{j(z,x_n(z))}{\|x_n\|_{H^1_0(Z)}^2}\hbox{d} z\\[.4pc]
&\geq \frac{c_4}{2}-\!\frac{\lambda_k}{2}\| y_n\|_2^2
-\frac{1}{2}\!\int_Z\gamma(z)y_n(z)^2\hbox{d} z
-\frac{\varepsilon}{2}\|y_n\|_2^2
-\frac{\xi_2|Z|}{\|x_n\|_{H^1_0(Z)}^2},
\end{align*}
for some $c_4>0$ (by the Poincar\'e inequality). Passing to the
limit as $n\rightarrow+\infty$, we obtain
\begin{equation*}
0 \geq \frac{1}{2}(c_4
-(\lambda_k+\|\gamma\|_{\infty}+\varepsilon)\|y\|_2^2)
\end{equation*}
and thus $y\not=0$.

Because of hypothesis H(j)(iv), we can find $M_3>0$ such that for
almost all $z\in Z$, all $|\zeta|\geq M_3$ and all $u\in\partial
j(z,\zeta)$, we have
\begin{equation*}
u\zeta-2j(z,\zeta) \leq -1.
\end{equation*}

On the other hand, as above, for almost all $z\in Z$ and all
$|\zeta|< M_3$, we have $|j(z,\zeta)|\leq \xi_3$ for some
$\xi_3>0$. Thus using hypothesis H(j)(iii), we see that for almost
all $z\in Z$, all $|\zeta|< M_3$ and all $u\in\partial
j(z,\zeta)$, we have
\begin{equation*}
|u\zeta-2j(z,\zeta)| \leq\xi_4,
\end{equation*}
for some $\xi_4>0$. Thus finally, for almost all $z\in Z$, all
$\zeta\in\mathbb{R}$ and all $u\in\partial j (z,\zeta)$, we have
\begin{equation} \label{eql_83}
u\zeta-2 j(z,\zeta) \leq \xi_4.
\end{equation}
Let $C\stackrel{d f}{=}\{z\in Z\hbox{\rm :}\ y(z)\not=0\}$.
Evidently $|C|_N>0$ (with $|\cdot|_N$ being the Lebesgue measure
on $\mathbb{R}^N$) and for all $z\in C$, we have that
\begin{equation} \label{eql11}
|x_n(z)|\longrightarrow+\infty \quad\hbox{as}\
n\rightarrow+\infty.
\end{equation}

From~(\ref{eql11}) and by virtue of Lemma~1 of \cite{TWb}, for a
given $\delta\in (0,|C|_N)$, we can find a measurable subset
$C_1\subseteq C$, such that $|C\!\setminus\!C_1|_N<\delta$ and
$|x_n(z)|\longrightarrow+\infty$ uniformly for all $z\in C_1$.
Using hypothesis H(j)(iv), we infer that
\begin{equation*}
\int_{C_1}(h_n(z)x_n(z)-2j(z,x_n(z)))\hbox{d} z
\longrightarrow-\infty \quad\hbox{as}\ n\rightarrow+\infty.
\end{equation*}
From~(\ref{eql_83}), we have
\begin{equation*}
h_n(z)x_n(z)-2j(z,x_n(z))\leq\xi_4 \quad\hbox{for a.a.}\ z\in
Z\!\setminus\!C_1,
\end{equation*}
and so
\begin{align*}
&\int_Z(h_n(z)x_n(z)-2j(z,x_n(z)))\hbox{d} z\\[.4pc]
&\quad\, = \int_{C_1}(h_n(z)x_n(z)-2j(z,x_n(z)))\hbox{d} z
+\xi_4|(Z\!\setminus\!C_1)_c|_N \longrightarrow-\infty,
\end{align*}
as $n\rightarrow+\infty$, which contradicts~(\ref{eq14}). This
proves the boundedness of $\{x_n\}\subseteq H^1_0(Z)$. So, passing
to a subsequence if necessary, we may assume that
\begin{align*}
x_n\stackrel{w}{\longrightarrow} x_0 \quad &\hbox{in}\ H^1_0(Z)\\[.4pc]
x_n\longrightarrow x_0               \quad &\hbox{in}\ L^2(Z),
\end{align*}
for some $x_0\in H^1_0(Z)$. From the choice of the sequence
$\{x_n\}_{n\geq 1}\subseteq H^1_0(Z)$, we have
\begin{align*}
&|\langle x_n^*,x_n-x_0\rangle_{H^1_0(Z)}|\\[.4pc]
&\quad\,=|\langle A(x_n),x_n-x_0\rangle_{H^1_0(Z)}
-\lambda_k\langle x_n^*,x_n-x\rangle_{H^1_0(Z)}\\[.4pc]
&\qquad\,-\int_Zh_n(z)(x_n(z)-x_0(z))\hbox{d} z|\\[.4pc]
&\quad\,\leq m_{\varphi}(x_n)\|x_n-x_0\|_{H^1_0(Z)}\longrightarrow
0 \quad\hbox{as} \ n\rightarrow+\infty.
\end{align*}
But
\begin{equation*}
\lambda_k\langle x_n^*,x_n-x\rangle_{H^1_0(Z)}
+\int_Zh_n(z)(x_n(z)-x_0(z))\hbox{d} z \longrightarrow
0\quad\hbox{as}\ n\rightarrow+\infty.
\end{equation*}
So it follows that
\begin{equation*}
\lim_{n\rightarrow+\infty}\langle A(x_n),x_n-x_0\rangle_{H^1_0(Z)} =  0.
\end{equation*}
From the maximal monotonicity of $A$, we have
\begin{equation*}
\langle A(x_n),x_n\rangle_{H^1_0(Z)} \longrightarrow \langle
A(x_0),x_0)\rangle_{H^1_0(Z)}
\end{equation*}
and so
\begin{equation*}
\|\nabla x_n\|_2 \longrightarrow \|\nabla x_0\|_2.
\end{equation*}

Because $\nabla x_n\stackrel{w}{\longrightarrow}\nabla x_0$ in
$L^2(Z;\mathbb{R}^N)$, from the Kadec--Klee property of Hilbert
spaces, we conclude that $\nabla x_n\longrightarrow\nabla x_0$ in
$L^2(Z;\mathbb{R}^N)$, hence $x_n\longrightarrow x_0$ in
$H^1_0(Z)$ as\break $n\rightarrow+\infty$.\hfill $\Box$
\end{proof}

Using this proposition, we can establish that $\bar{\varphi}$
satisfies the nonsmooth C-condition.

\begin{proposition} \label{prop6}$\left.\right.$\vspace{.5pc}

\noindent If hypotheses {\rm H(j)} hold{\rm ,} then
$\bar{\varphi}$ satisfies the nonsmooth ${\rm C}$-condition.
\end{proposition}

\begin{proof}
Let $c\in \mathbb{R}$ and let $\{u_n\}_{n\geq 1}\subseteq
\bar{H}_0$ be a sequence, such that
\begin{equation*}
\bar{\varphi}(u_n)\longrightarrow c \quad\hbox{and}\quad
(1+\|u_n\|_{H^1_0(Z)})m_{\bar{\varphi}}(u_n)\longrightarrow 0
\quad\hbox{as}\ n\rightarrow+\infty.
\end{equation*}
As before we can find $\bar{v}_n^*\in\partial\bar{\varphi}(u_n)$,
such that $m_{\bar{\varphi}}(u_n)=\|\bar{v}_n^*\|_{H^{-1}(Z)}$. By
virtue of~(\ref{eq10}), we can find
$v_n^*\in\partial\varphi(u_n+\vartheta(u_n))$, such that
\begin{equation*}
\bar{v}_n^* = p_{\bar{H}_0^*}v_n^* \quad\forall n\geq 1.
\end{equation*}
Recall that, by Proposition~\ref{prop3}, $0\in
p_{\hat{H}^*}\partial\varphi(u_n+\vartheta(u_n))$, for $n\geq 1$.
Then using hypothesis H(j)(v) we have
$m_{\varphi}(u_n+\vartheta(u_n))\leq m_{\bar{\varphi}}(u_n).$
Therefore,
\begin{equation*}
(1+\|u_n\|_{H^1_0(Z)})\|v_n^*\|_{H^1_0(Z)} \longrightarrow 0
\quad\hbox{with}\ v_n^*\in\partial\varphi(u_n+\vartheta(u_n)).
\end{equation*}

But from Proposition~\ref{prop5}, we know that $\varphi$ satisfies
the nonsmooth C-condition. So we can extract a subsequence of
$\{u_n\}_{n\geq 1}$, which is strongly convergent. This proves the
proposition.\hfill $\Box$
\end{proof}

\begin{proposition} \label{prop7}$\left.\right.$\vspace{.5pc}

\noindent If hypotheses ${\rm H(j)}$ hold{\rm ,} then $\psi$ is
bounded below.
\end{proposition}

\begin{proof}
We show that $-(\varphi|_{\bar{H}_0})$ is bounded below. Then
because $-\psi=\bar{\varphi}\leq(\varphi|_{\bar{H}_0})$, we can
conclude that $\psi$ is bounded below. To this end we proceed by
contradiction. Suppose that $-(\varphi|_{\bar{H}_0})$ is not
bounded below. Then we can find $x_n\in\bar{H}_0$, such that
\begin{equation} \label{eql_80}
\varphi(x_n) \geq n \quad\forall n\geq 1
\end{equation}
and
\begin{equation*}
\|x_n\|_{H^1_0(Z)} \longrightarrow +\infty.
\end{equation*}
By virtue of hypothesis H(j)(vii), for a given $\varepsilon>0$, we
can find $M_4=M_4(\varepsilon)>0$, such that for almost all $z\in
Z$ and all $|\zeta|\geq M_4$, we have
\begin{equation*}
-\frac{\varepsilon}{2}\zeta^2 \leq j(x,\zeta).
\end{equation*}
On the other hand, as before via Lebourg mean value theorem, we
can find $\xi_5>0$, such that for almost all $z\in Z$ and all
$|\zeta|\leq M_4$, we have
\begin{equation*}
|j(z,\zeta)| \leq \xi_5.
\end{equation*}
So finally we see that
\begin{equation} \label{eq15}
-\frac{\varepsilon}{2}\zeta^2-\xi_5 \leq j(z,\zeta)
\quad\hbox{for a.a.}\ z\in Z,\ \hbox{all}\ \zeta\in\mathbb{R}.
\end{equation}

Let
\begin{equation*}
x_n=\bar{x}_n+\bar{\hskip -.05pc \bar{x}}_n, \quad\hbox{with}\
\bar{x}_n\in\bar{H}, \bar{\hskip -.05pc \bar{x}}_n\in
E(\lambda_k), \quad n\geq 1.
\end{equation*}
First assume that
\begin{equation} \label{eql321}
\frac{\|\nabla\bar{x}_n\|_2}{\|\nabla x_n\|_2}
 \longrightarrow \mu\not=0
\quad\hbox{as}\ n\rightarrow+\infty.
\end{equation}
Exploiting the orthogonality relations, the fact that $\|\nabla
\bar{\hskip -.05pc \bar{x}}_n\|_2^2 =\lambda_k\|\bar{\hskip -.05pc
\bar{x}}_n\|_2^2$ and estimate~(\ref{eq15}), we have
\begin{align*}
\varphi(x_n) &=\frac{1}{2}\|\nabla \bar{x}_n\|_2^2
-\frac{\lambda_k}{2}\|\bar{x}_n\|_2^2
-\int_Zj(z,x_n(z))\hbox{d} z\\[.4pc]
&\leq\frac{1}{2}\|\nabla\bar{x}_n\|_2^2
-\frac{\lambda_k}{2}\|\bar{x}_n\|_2^2
+\frac{\varepsilon}{2}\|x_n\|_2^2 +\xi_5|Z|_N.
\end{align*}
Thus from the variational characterization of the eigenvalues we
get
\begin{align} \label{eq16}
\varphi(x_n) &\leq
\frac{1}{2}\left(1-\frac{\lambda_k}{\lambda_{k-1}}\right)
\|\nabla\bar{x}_n\|_2^2 +\frac{\varepsilon}{2}\|x_n\|_2^2
+\xi_5|Z|_N\nonumber\\[.4pc]
&\leq \frac{1}{2}\|\nabla x_n\|_2^2 \left(
\left(1-\frac{\lambda_k}{\lambda_{k-1}}\right)
\frac{\|\nabla\bar{x}_n\|_2^2}{\|\nabla x_n\|_2^2}
+\frac{\varepsilon}{\lambda_k} \right) +\xi_5|Z|_N.
\end{align}

Since by hypothesis, we have that $\|\nabla
x_n\|_2\longrightarrow+\infty$, so from~(\ref{eql321}) and
recalling that $\lambda_{k-1}<\lambda_k$, by passing to the limit
as $n\rightarrow+\infty$ in~(\ref{eq16}), we see that
$\varphi(x_n)\longrightarrow-\infty$ as $n\rightarrow+\infty$, a
contradiction to~(\ref{eql_80}).

Next assume that
\begin{equation} \label{eql322}
\frac{\|\nabla\bar{x}_n\|_2}{\|\nabla x_n\|_2} \longrightarrow 0
\quad\hbox{as}\ n\rightarrow+\infty.
\end{equation}
By virtue of hypothesis H(j)(iv), for a given $\eta>0$, we can
find $M_5=M_5(\eta)>0$, such that for almost all $z\in Z$, all
$|\zeta|\geq M_5$ and all $u\in\partial j(z,\zeta)$, we have
\begin{equation} \label{eq17}
u\zeta-2j(z,\zeta) \leq -\eta.
\end{equation}

From p.~48 of \cite{Cl}, we know that for almost all $z\in Z$ and
all $\zeta>0$, the function
$\zeta\longmapsto\frac{j(z,\zeta)}{\zeta^2}$ is locally Lipschitz
and we have that
\begin{align*}
\partial\left(\frac{j(z,\zeta)}{\zeta^2}\right)
&= \frac{\zeta^2\partial j(z,\zeta)-2\zeta
j(z,\zeta)}{\zeta^4}\\[.4pc]
&= \frac{\zeta\partial j(z,\zeta)-2j(z,\zeta)}{\zeta^3}.
\end{align*}
Using~(\ref{eq17}), we see that for almost all $z\in Z$, all
$\zeta\geq M_5$ and all
$v\in\partial\!\left(\frac{j(z,\zeta)}{\zeta^2}\right)$, we have
\begin{equation*}
v \leq -\frac{\eta}{\zeta^3}.
\end{equation*}
Since for all $z\in Z\!\setminus\!E$, with $|E|_N=0$, the function
$\zeta\longmapsto\frac{j(z,\zeta)}{\zeta^2}$ is locally Lipschitz
on $[M_5,+\infty)$, it is differentiable at every
$x\in[M_5,+\infty)\!\setminus\!L(z)$ with some $|L(z)|_1=0$. We
set
\begin{equation*}
\xi_0(z,\zeta)\stackrel{d f}{=} \begin{cases} \displaystyle
\frac{{\rm d}}{{\rm d} z}\left(\frac{j(z,\zeta)}{\zeta^2}\right),
&\hbox{if} \quad x\in[M_5,+\infty)\!\setminus\!L(z),\\[.4pc]
0, &\hbox{if} \quad x\in L(z).
\end{cases}
\end{equation*}
For all $z\in Z\!\setminus\!E$ and all
$\zeta\in[M_5,+\infty)\!\setminus\!L(z)$, we have
$\xi_0(z,\zeta)\in\partial\!\left(\frac{j(z,\zeta)}{\zeta^2}\right)$
and so
\begin{equation*}
\xi_0(z,\zeta) \leq -\frac{\eta}{\zeta^3}.
\end{equation*}
Let $\tau,\bar{\tau}\in [M_5,+\infty)$, with $\tau\leq\bar{\tau}$.
Integrating the above inequality with respect to
$\zeta\in\mathbb{R}$ on the interval $[\tau,\bar{\tau}]$, we
obtain
\begin{equation*}
\frac{j(z,\bar{\tau})}{\bar{\tau}^2} -\frac{j(z,\tau)}{\tau^2}
\leq -\frac{\eta}{2}
\left(\frac{1}{\tau^2}-\frac{1}{\bar{\tau}^2}\right).
\end{equation*}
Let $\bar{\tau}\longrightarrow+\infty$ and using hypothesis
H(j)(vii), we obtain
\begin{equation*}
\frac{j(z,\tau)}{\tau^2} \geq \frac{\eta}{2\tau^2}
\end{equation*}
and thus, for almost all $z\in Z$ and all $\tau\geq M_5$, we have
that
\begin{equation*}
j(z,\tau) \geq \frac{\eta}{2}.
\end{equation*}
Since $\eta>0$ was arbitrary, it follows that
$j(z,\tau)\longrightarrow+\infty$ as $\tau\rightarrow+\infty$
uniformly for almost all $z\in Z$. Similarly we show that
$j(z,\tau)\longrightarrow+\infty$, $\tau\rightarrow-\infty$
uniformly for almost all $z\in Z$. Therefore
\begin{equation} \label{eql_214}
j(z,\zeta)\longrightarrow+\infty \quad \hbox{as}\
|\zeta|\rightarrow+\infty \ \hbox{uniformly for a.a.}\ z\in Z.
\end{equation}

From \cite{BBF} (see the proof of Lemma~3.2) or \cite{ST}, we have
that for a given $\delta>0$ we can find $\xi_{\delta}>0$, such
that
\begin{equation*}
|\{ z\in Z\hbox{\rm :}\ |v(z)|<\xi_{\delta}\|v\|_{H^1_0(Z)}\} |_N
< \delta \quad\forall v\in E(\lambda_k).
\end{equation*}
Let us set
\begin{equation*}
C_n\stackrel{d f}{=} \{ z\in Z\hbox{\rm :}\ |\bar{\hskip -.05pc
\bar{x}}_n(z)| \geq\xi_{\delta}\|\bar{\hskip -.05pc
\bar{x}}_n\|_{H^1_0(Z)}\} \quad\forall n\geq 1.
\end{equation*}
We have that $|Z\!\setminus\!C_n|_N<\delta$.

Let us establish that the sequence $\{\bar{x}_n\}_{n\geq 1}$ is
bounded in $H^1_0(Z)$. Arguing by contradiction, assume that along
a relabeled subsequence we have
$\|\bar{x}_n\|_{H^1_0(Z)}\rightarrow +\infty$. Recall that
\begin{equation*}
\varphi(x_n) \leq \frac{1}{2}
\left(1-\frac{\lambda_k}{\lambda_{k-1}}\right)
\|\nabla\bar{x}_n\|_2^2 -\int_Zj(z,x_n(z))\hbox{d} z \quad\forall
n\geq 1.
\end{equation*}
Using~(\ref{eql_214}) we find a constant $\alpha>0$ such that
$j(z,\zeta)\geq 0$ a.a. $z\in Z$, $\forall |\zeta|>\alpha$. We see
that
\begin{align*}
-\int_Zj(z,x_n(z))\hbox{d} z &=
-\int_{\{z\in Z\hbox{\rm :}\ |x_n(z)|>\alpha\}}j(z,x_n(z))\hbox{d} z\\[.4pc]
&\quad\,
-\int_{\{z\in Z\hbox{\rm :}\ |x_n(z)|\leq \alpha\}}j(z,x_n(z))\hbox{d} z\\[.4pc]
&\leq
-\int_{\{z\in Z\hbox{\rm :}\ |x_n(z)|\leq \alpha\}}j(z,x_n(z))\hbox{d} z\\[.4pc]
&\leq k_0|Z|_N,
\end{align*}
where $k_0>0$ is a constant. In writing the last inequality above
we made use of assumption H(j)(iii) and Lebourg's mean value
theorem. Then we obtain the estimate
\begin{equation*}
\varphi(x_n)\leq
\frac{1}{2}\left(1-\frac{\lambda_k}{\lambda_{k-1}}\right)
\|\nabla\bar{x}_n\|^2_2+k_0|Z|_N.
\end{equation*}
Since $\lambda_{k-1}<\lambda_k$ and we supposed that
$\|\bar{x}_n\|_{H^1_0(Z)}\rightarrow +\infty$, we arrive at the
conclusion that $\varphi(x_n)\rightarrow -\infty$ as $n\rightarrow
+\infty$. This contradicts relation~(\ref{eql_80}), and thus the
sequence $\{\bar{x}_n\}_{n\geq 1}$ is bounded in $H^1_0(Z)$.

Because $\bar{x}_n\in \bar{H}$ and the space $\bar{H}$ is finite
dimensional, it follows from the boundedness of the sequence
$\{\bar{x}_n\}_{n\geq 1}\subseteq H^1_0(Z)$ that we can find
$c_5>0$ such that
\begin{equation*}
|\bar{x}_n(z)|\leq c_5 \quad\forall z\in Z, n\geq 1.
\end{equation*}
From~(\ref{eql_214}), we know that for a given $\eta_1>0$ we can
find $M_6=M_6(\eta_1)>0$ such that
\begin{equation*}
j(z,\zeta)\geq\eta_1 \quad\hbox{for a.a.}\ z\in Z, \ \hbox{all}\
|\zeta|\geq M_6.
\end{equation*}
Let
\begin{equation*}
D_n\stackrel{d f}{=} \{z\in Z\hbox{\rm :}\ |x_n(z)|\geq M_6\}
\quad\forall n\geq 1.
\end{equation*}
If $z_0\in C_n$, then
\begin{equation*}
|x_n(z_0)| \geq |\bar{\hskip -.05pc \bar{x}}_n(z_0)|
-|\bar{x}_n(z_0)| \geq \xi_{\delta}\|\bar{\hskip -.05pc
\bar{x}}_n\|_{H^1_0(Z)}-c_5.
\end{equation*}
Because of~(\ref{eql322}), we must have that $\|\bar{\hskip -.05pc
\bar{x}}_n\|_{H^1_0(Z)}\longrightarrow+\infty$ as
$n\rightarrow+\infty$. So there exists $n_0\geq 1$ large enough
such that
\begin{equation*}
\xi_{\delta}\|\bar{\hskip -.05pc \bar{x}}_n\|_{H^1_0(Z)} -c_5 \geq
M_6 \quad \forall n\geq n_0
\end{equation*}
and so $z_0\in D_n$ for $n\geq n_0$, i.e.
\begin{equation*}
C_n\subseteq D_n \quad\forall n\geq n_0.
\end{equation*}
Then using (\ref{eq15}) and the fact that $|Z\!\setminus\!
D_n|<\delta$ (as $Z\!\setminus\!D_n\subseteq Z\!\setminus\!C_n$),
for $n\geq n_0$, we have that
\begin{align*}
\int_Zj(z,x_n(z))\hbox{d} z &= \int_{D_n}j(z,x_n(z))\hbox{d} z
+\int_{Z\setminus D_n}j(z,x_n(z))\hbox{d} z\\[.4pc]
&\geq \eta_1|D_n|_N -\left(\frac{\varepsilon}{2}M_6^2+\xi_5\right)
|Z\!\setminus\!D_n|_N\\[.4pc]
&\geq \eta_1|D_n|_N -\left(\frac{\varepsilon}{2}M_6^2+\xi_5\right)
\delta\\[.4pc]
&\geq \eta_1(|Z|_N-\delta)
-\left(\frac{\varepsilon}{2}M_6^2+\xi_5\right) \delta,
\end{align*}
so
\begin{equation*}
\liminf_{n\rightarrow+\infty}\int_Zj(z,x_n(z))\hbox{d} z \geq
\eta_1(|Z|_N-\delta)
-\left(\frac{\varepsilon}{2}M_6^2+\xi_5\right) \delta.
\end{equation*}
Since $\delta>0$ is arbitrary, we let $\delta\searrow 0$. We
obtain
\begin{equation*}
\liminf_{n\rightarrow+\infty}\int_Zj(z,x_n(z))\hbox{d} z \geq
\eta_1|Z|_N.
\end{equation*}
Because $\eta_1>0$ is arbitrary, we conclude that
\begin{equation*}
\int_Zj(z,x_n(z))\hbox{d} z \longrightarrow+\infty \quad\hbox{as}\
n\rightarrow+\infty.
\end{equation*}
From the choice of the sequence $\{x_n\}_{n\geq
1}\subseteq\bar{H}_0$, we have
\begin{align*}
n \leq \varphi(x_n) \leq
\frac{1}{2}\left(1-\frac{\lambda_k}{\lambda_{k-1}}\right)
\|\nabla\bar{x}_n\|^2_2 -\int_Z j(z,x_n(z))\hbox{d} z
\longrightarrow -\infty,
\end{align*}
a contradiction. Therefore $-(\varphi|_{\bar{H}_0})$ is bounded
below and so $\psi$ is bounded below too.\quad \hfill $\Box$
\end{proof}

Now we are ready for our multiplicity result.

\begin{theorem}[\!] \label{thm8}
If hypotheses ${\rm H(j)}$ hold{\rm ,} then problem~${\rm (HVI)}$
has at least two nontrivial solutions.
\end{theorem}

\begin{proof}
If $\inf_{\bar{H}_0}\psi=0$ (remark that $\psi(0)=0$), then by
virtue of Proposition~\ref{prop4}, all $x\in V$ with
$\|x\|_V\leq\delta$ are critical points of $\psi$.

So assume that $\inf_{\bar{H}_0}\psi<0$. By virtue of
Propositions~\ref{prop4}, \ref{prop6} and~\ref{prop7} and since
$\psi(0)=0$, we can apply Theorem~\ref{thm1} and obtain two
nontrivial critical points $x_1,x_2\in\bar{H}_0$ of $\psi$, i.e.
\begin{equation*}
0\in\partial\psi(x_i) \quad\hbox{for}\ i=1,2.
\end{equation*}
Hence $0\in\partial\bar{\varphi}(x_i)$ for $i=1,2$ and finally
from~(\ref{eq10}), we have
\begin{equation*}
0\in p_{\bar{H}_0^*}\partial\varphi(x_i+\vartheta(x_i))
\quad\hbox{for}\ i=1,2
\end{equation*}
(see the definition of $\vartheta$ in Proposition~\ref{prop3}).
Recall that $0\in
p_{\hat{H}^*}\partial\varphi(x_i+\vartheta(x_i))$ for $i=1,2$,
where $\hat{H}=\bar{H}_0^{\perp}$. Therefore we have
$0\in\partial\varphi(x_i+\vartheta(x_i))$ for $i=1,2$. Set
$u_i=x_i+\vartheta(x_i)$ for $i=1,2$. Then $u_1$ and $u_2$ are two
nontrivial critical points of $\varphi$, and thus they are two
nontrivial solutions of~(HVI).\hfill $\Box$
\end{proof}

\begin{remark}
{\rm An example of a nonsmooth locally Lipschitz potential
satisfying hypotheses H(j) is the following. For simplicity we
drop the $z$-dependence (see figure~1 for $j$ and figure~2 for the
Clarke subdifferential $\partial j$).}
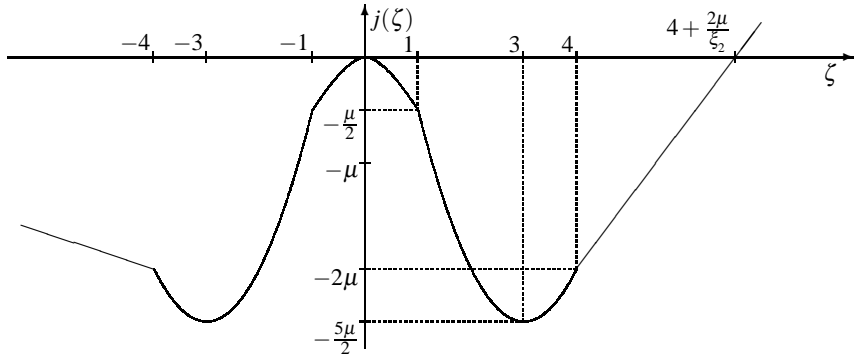
\begin{figure}[h]\vspace{1pc}
\begin{center}
\begin{picture}(330,130)(0,0)
\put(5,110){\vector(1,0){320}}   
\put(140,0){\vector(0,1){130}}   
\put(60,30){\line(-3,1){50}}
\bezier{300}(80,10)(70,10)(60,30)   
\bezier{450}(120,90)(100,10)(80,10)   
\bezier{150}(140,110)(133,110)(120,90)   
\bezier{150}(140,110)(147,110)(160,90)   
\bezier{450}(160,90)(180,10)(200,10)   
\bezier{300}(200,10)(210,10)(220,30)   
\put(220,30){\line(3,4){70}}
\put(138,90){\line(1,0){4}}
\put(138,90){\makebox(0,0)[rt]{{\scriptsize $-\frac{\mu}{2}$}}}
\multiput(143,90)(2,0){9}{\line(1,0){1}}
\put(138,70){\line(1,0){4}}
\put(138,70){\makebox(0,0)[rt]{{\scriptsize $-\mu$}}}
\put(138,30){\line(1,0){4}}
\put(138,30){\makebox(0,0)[rt]{{\scriptsize $-2\mu$}}}
\multiput(143,30)(2,0){39}{\line(1,0){1}}
\put(138,10){\line(1,0){4}}
\put(138,10){\makebox(0,0)[rt]{{\scriptsize $-\frac{5\mu}{2}$}}}
\multiput(143,10)(2,0){28}{\line(1,0){1}}
\put(142,129){\makebox(0,0)[lt]{{\scriptsize $j(\zeta)$}}}
\put(60,108){\line(0,1){4}}
\put(59,112){\makebox(0,0)[rb]{{\scriptsize $-4$}}}
\put(80,108){\line(0,1){4}}
\put(79,112){\makebox(0,0)[rb]{{\scriptsize $-3$}}}
\put(120,108){\line(0,1){4}}
\put(119,112){\makebox(0,0)[rb]{{\scriptsize $-1$}}}
\put(160,108){\line(0,1){4}}
\put(159,112){\makebox(0,0)[rb]{{\scriptsize $1$}}}
\multiput(160,90)(0,2){9}{\line(0,1){1}}
\put(200,108){\line(0,1){4}}
\put(199,112){\makebox(0,0)[rb]{{\scriptsize $3$}}}
\multiput(200,10)(0,2){49}{\line(0,1){1}}
\put(220,108){\line(0,1){4}}
\put(219,112){\makebox(0,0)[rb]{{\scriptsize $4$}}}
\put(280,108){\line(0,1){4}}
\put(279,112){\makebox(0,0)[rb]{{\scriptsize $4+\frac{2\mu}{\xi_{{}_{\mbox{\tiny 2}}}}$}}}
\multiput(220,30)(0,2){39}{\line(0,1){1}}
\put(319,108){\makebox(0,0)[rt]{{\scriptsize $\zeta$}}}
\end{picture}
\end{center}\vspace{-1.3pc}
\caption{Potential function.}
\end{figure}
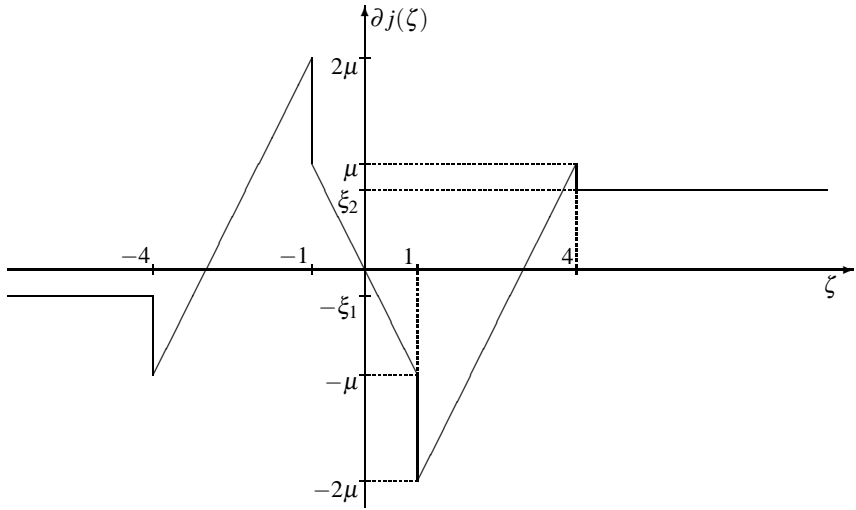
\begin{figure}[h]
\begin{center}
\begin{picture}(330,190)(0,0)
\put(5,90){\vector(1,0){320}}   
\put(140,0){\vector(0,1){190}}   
\put(5,80){\line(1,0){55}} 
\put(60,80){\line(0,-1){30}} 
\put(60,50){\line(1,2){60}} 
\put(120,170){\line(0,-1){40}} 
\put(120,130){\line(1,-2){40}} 
\put(160,50){\line(0,-1){40}} 
\put(160,10){\line(1,2){60}} 
\put(220,130){\line(0,-1){10}} 
\put(220,120){\line(1,0){95}} 
\put(138,170){\line(1,0){4}}
\put(138,170){\makebox(0,0)[rt]{{\scriptsize $2\mu$}}}
\put(138,130){\line(1,0){4}}
\put(138,130){\makebox(0,0)[rt]{{\scriptsize $\mu$}}}
\multiput(143,130)(2,0){39}{\line(1,0){1}}
\put(138,120){\line(1,0){4}}
\put(138,120){\makebox(0,0)[rt]{{\scriptsize $\xi_2$}}}
\multiput(143,120)(2,0){39}{\line(1,0){1}}
\put(138,80){\line(1,0){4}}
\put(138,80){\makebox(0,0)[rt]{{\scriptsize $-\xi_1$}}}
\put(138,50){\line(1,0){4}}
\put(138,50){\makebox(0,0)[rt]{{\scriptsize $-\mu$}}}
\multiput(143,50)(2,0){9}{\line(1,0){1}}
\put(138,10){\line(1,0){4}}
\put(138,10){\makebox(0,0)[rt]{{\scriptsize $-2\mu$}}}
\multiput(143,10)(2,0){9}{\line(1,0){1}}
\put(142,189){\makebox(0,0)[lt]{{\scriptsize $\partial j(\zeta)$}}}
\put(60,88){\line(0,1){4}}
\put(59,92){\makebox(0,0)[rb]{{\scriptsize $-4$}}}
\put(120,88){\line(0,1){4}}
\put(119,92){\makebox(0,0)[rb]{{\scriptsize $-1$}}}
\put(160,88){\line(0,1){4}}
\put(159,92){\makebox(0,0)[rb]{{\scriptsize $1$}}}
\multiput(160,50)(0,2){19}{\line(0,1){1}}
\put(220,88){\line(0,1){4}}
\put(219,92){\makebox(0,0)[rb]{{\scriptsize $4$}}}
\multiput(220,93)(0,2){14}{\line(0,1){1}}
\put(319,88){\makebox(0,0)[rt]{{\scriptsize $\zeta$}}}
\end{picture}
\end{center}\vspace{-1.3pc}
\caption{Subdifferential of the potential.}\vspace{1pc}
\end{figure}
{\rm \begin{equation*} j(\zeta)\stackrel{d f}{=}\begin{cases}
-\xi_1\zeta-2\mu-4\xi_1, &\hbox{if} \quad \zeta <-4\\[.6pc]
\displaystyle\frac{\mu}{2}\zeta^2+3\mu\zeta+2\mu, &\hbox{if} \quad -4\leq\zeta<-1\\[.6pc]
\displaystyle-\frac{\mu}{2}\zeta^2, &\hbox{if} \quad -1\leq\zeta<1\\[.6pc]
\displaystyle\frac{\mu}{2}\zeta^2-3\mu\zeta+2\mu, &\hbox{if} \quad 1\leq \zeta<4\\[.6pc]
\xi_2\zeta-2\mu-4\xi_2, &\hbox{if} \quad 4\leq \zeta.
\end{cases}
\end{equation*}
Here $\lambda_k-\lambda_m<\mu< \min\{\lambda_k-\lambda_{m-1},\
\lambda_{k+1}-\lambda_k\} $, $0<\xi_1,\xi_2<\mu$. All the
\hbox{assumptions~(i)--(vii)} in H(j) are verified. For instance,
assumption H(j)(v) holds with $l(z)\equiv\mu$. In this case we
have resonance at $\pm\infty$ since
$\frac{j(\zeta)}{\zeta^2}\longrightarrow 0$ as
$|\zeta|\rightarrow+\infty$. Another possibility is the function
$j(x)=\max{\big\{\frac{\xi}{2}x^2+c|x|,\frac{\xi}{2}|x|\big\}}$
with $\xi<\lambda_{k+1}-\lambda_{k},\;c\leq \frac{\xi}{2}.$}
\end{remark}

\section*{Acknowledgement}

The authors wish to thank the referee for his constructive
remarks. The first author (LG) is an award holder of the NATO
Science Fellowship Programme, which was spent in the National
Technical University of Athens.

\end{document}